\documentclass[12pt]{article}

\usepackage[cp1251]{inputenc}
\usepackage[english]{babel}
\usepackage{amsmath}
\usepackage{amsthm}
\usepackage{amsxtra}
\usepackage{amsfonts,amssymb}

\newtheorem{Th}{Theorem}
\newtheorem{Prop}[Th]{Proposition}
\newtheorem{Lm}[Th]{Lemma}

\theoremstyle{definition}

\date{}
\begin{document}

\title{
    \bf
Asymptotics of Plancherel measures for $GL(n,q)$.}

\author{A.\,Dudko}

\date{}
\maketitle

\begin{abstract}
We introduce the Plancherel measures on the sets of partition collections
$\{\Lambda_\phi\}$, which parameterize irreducible representations of
$GL(n,q)$. Consider the Plancherel random collection $\{\Lambda_\phi\}$. We
prove that as $n\rightarrow\infty$, the random partitions $\Lambda_\phi$
converge in finite dimensional distribution to independent random
variables. We give explicit formulas for the corresponding limit
distributions.
\end{abstract}

\paragraph{1. Introduction.} Let $GL(n,q)$ be the group of all invertible $n\times n$
matrices over $\mathbb{F}_q$, where $q$ is a prime power. Denote by
$\Phi$ the set of all irreducible polynomials over $\mathbb{F}_q$
with unit highest coefficient and nonzero constant term. The
irreducible representations of $GL(n,q)$ are parameterized by
collections $\Lambda=\{\Lambda_\phi\}$ of partitions, where $\phi\in
\Phi$ and
$$ |\Lambda|=\sum\limits_\phi\left|\Lambda_\phi\right|deg\,\phi=n.$$
 Denote by $\mathbb{L}_n$ the set of collections
 $\Lambda=\{\Lambda_\phi\}$ with $|\Lambda|=n$.
 Define the
 Plancherel measure $\mu_n$ on $\mathbb{L}_n$ as follows:
 \begin{eqnarray}
 \label{Plansh}\mu_n(\Lambda)=\frac{d^2_\Lambda}{|GL(n,q)|},\end{eqnarray}
  where $d_\Lambda$ is the degree of the irreducible representation
  of $GL(n,q)$ corresponding to $\Lambda\in \mathbb{L}_n$.
The Plancherel measure defines relative dimensions of regular
representation isotopic components.

  For $0<v<q$ and $q>1$ denote the measure $ M_{v,q}$ on the set of all partitions of all
numbers as follows:
\begin{eqnarray}\label{nu}
  M_{v,q}(\lambda)=\prod\limits_{r=1}^\infty\left(1-vq^{-r}\right)^r
 q^{2n(\lambda')+|\lambda|}\prod\limits_{x\in \lambda}
 \left(q^{h(x)}-1\right)^{-2}v^{|\lambda|},
 \end{eqnarray} where $h(x)$ is the
 hook-length and
 $n(\lambda')=\sum\limits_i\lambda_i(\lambda_i-1)/2$.
We show, that $ M_{v,q}$ is probabilistic. As we shall see, the
measure $ M_{v,q}$ is a specialization of the Schur measure (see
\cite{Ok}).

  Let $\Lambda$ be a random value in
  $\mathbb{L}_n$ distributed by (\ref{Plansh}). Then each
   $\phi\in\Phi$ denotes the random value $\Lambda_\phi$ in the set of partitions.
    We consider asymptotical behavior of $\Lambda_\phi$ with increase of
    $n$. The main result of this paper is the following
   \begin{Th}\label{main}
   Let $\lambda_1,\lambda_2,\ldots, \lambda_k$ be arbitrary
   partitions. Let $\phi_1,\phi_2,\ldots,\phi_k\in \Phi$ be different
   polynomials. Then there exists the limit:
\begin{eqnarray}
\lim_{n\rightarrow \infty}\mu_n(\{\Lambda\in\mathbb{L}_n
:\Lambda_{\phi_j}=\lambda_j,j=1,\ldots,k\})= \prod\limits_{j=1}^k
M_{1,\,q^{deg\,\phi_j}}(\lambda_j).
\end{eqnarray}
   \end{Th}
The proof is based on the cycle index technic (see \cite{F1},
\cite{F2}). Asymptotical behavior of Plancherel measures for
symmetric groups is well studied. See for example \cite{Bor-Ok-Ol}.

I am grateful to N.I. Nessonov for the statement of the problem and
constant attention to my work. I also would like to thank A. Borodin for
the important remarks on the result of this paper. After this paper was
posted to arxiv, the author knew from J. Fulman, that the results of this
paper with proofs are contained in \cite{F3}.

\paragraph{2. Notions.}
Make some preliminaries. Let
$\lambda=(\lambda_1,\lambda_2,\ldots,\lambda_k)$ be arbitrary
partition. By definition, $n(\lambda)=\sum\limits
\lambda'_j(\lambda'_j-1)/2$. The hook-length in $x=(i,j)\in \lambda$
is equal to
\begin{eqnarray}h(x)=\lambda_i+\lambda'_j-i-j+1.
\end{eqnarray}
Since $n(\lambda)=\sum\limits_{(i,j)\in \lambda}(\lambda'_j-i)$,
then
\begin{eqnarray}\label{sum of h}
\sum\limits_{x\in \lambda}h(x)=n(\lambda)+n(\lambda')+|\lambda|.
\end{eqnarray}
From \cite{Mac}, I.3, ex. 2, using (\ref{sum of h}), one can deduce the following specialization
 formula for the Schur functions:
 \begin{eqnarray}\label{S spec}
 s_\lambda(q^{-1},q^{-2},\ldots)=q^{n(\lambda')}
 \prod\limits_{x\in \lambda}
 \left(q^{h(x)}-1\right)^{-1}\;\; (\text{see also }\cite{Mac}, \text{ IV.6}).
 \end{eqnarray}
By definition of the Schur measure
$\mathfrak{M}$,\begin{eqnarray}\label{S measure}
\mathfrak{M}(\lambda)=\prod\limits_{i,j}(1-x_iy_j)s_\lambda(x)s_\lambda(y).
\end{eqnarray} From (\ref{nu}), (\ref{S spec}) and (\ref{S measure})
the next statement follows.
\begin{Prop} Let $0<v<q$ and $q>1$. If we set $x_i=q^{-i},y_j=vq^{1-j}$
then the Schur measure $\mathfrak{M}$ specializes to $ M_{v,q}$.
\end{Prop}
 By the Cauchy identity for the Schur functions, it follows that
 $M_{v,q}$ is probabilistic.  Recall now the formula for the degree of the
irreducible representation of $GL(n,\mathbb{F}_q)$ corresponding to
$\Lambda=\{\Lambda_\phi\}$:
\begin{eqnarray}\label{d lamda}d_\Lambda=\prod\limits_{i=1}^n(q^i-1)\prod\limits_\phi
 q^{n(\lambda')deg\,\phi}\prod\limits_{x\in \lambda}
 \left(q^{h(x)deg\,\phi}-1\right)^{-1}\;\;(\text{see }\cite{Mac}, \text{IV}.6).\end{eqnarray}
By (\ref{Plansh}), (\ref{S spec}) and (\ref{d lamda}), one has:
\begin{eqnarray}\mu_n(\Lambda)=q^{-n(n-1)/2}\prod\limits_{i=1}^n(q^i-1)
\prod\limits_\phi
s^2_{\Lambda_\phi}(q^{-deg\,\phi},q^{-2deg\,\phi},\ldots).
\end{eqnarray}

\paragraph{3. The proof of the main result.}
 Let $0<v<1$ and $q>1$.
 Introduce the measure $P_{v,q}$ on $\mathbb{L}=\bigcup\limits_{n=0}^\infty
  \mathbb{L}_n$ as follows:
  \begin{eqnarray}\label{M u}
    P_{v,\,q}(\Lambda)=\prod\limits_{r=0}^\infty\left(1-vq^{-r}\right)
    q^{n(n+1)/2}\prod\limits_{i=1}^n(q^i-1)^{-1}v^n\mu_n(\Lambda)\text{
    for }\Lambda\in \mathbb{L}_n.
  \end{eqnarray}
 It follows from the Euler's identity (see \cite{A}, page 19), that
 \begin{eqnarray}\label{eiler}
    \sum\limits_{n=0}^\infty
    q^{n(n+1)/2}\prod\limits_{i=1}^n(q^i-1)^{-1}v^n=
    \prod\limits_{r=0}^\infty\left(1-vq^{-r}\right)^{-1}.
  \end{eqnarray}
  Therefore, $P_{v,\,q}$ is probabilistic.
For each $A\subset\mathbb{L}$
  one has:
\begin{eqnarray}\label{M u 2}
    P_{v,\,q}(A)\prod\limits_{r=0}^\infty\left(1-vq^{-r}\right)^{-1}=
    \sum\limits_{n=0}^\infty q^{n(n+1)/2}\prod\limits_{i=1}^n(q^i-1)^{-1}v^n\mu_n
    (A\cap \mathbb{L}_n).
  \end{eqnarray} The right side of (\ref{M u 2}) is the Taylor
  expansion for the left side in $|v|<1$.
    Consequently, for each $n$,
\begin{eqnarray}\label{M u 3}
   \prod\limits_{i=1}^n(1-q^{-i})^{-1}\mu_n(A\cap \mathbb{L}_n)=
   [v^n]
\prod\limits_{r=0}^\infty\left(1-vq^{-r}\right)^{-1}P_{v,\,q}(A),
  \end{eqnarray} where $[v^n]f(v)$ stands for the
  coefficient of the $n$-th power in the Taylor expansion for $f$.
Further, by (\ref{Plansh}), (\ref{nu}), (\ref{d lamda}) and (\ref{M
u}), for each $\Lambda\in \mathbb{L}_n$,
 \begin{eqnarray}\label{M u 4}
\prod\limits_{r=0}^\infty\left(1-vq^{-r}\right)^{-1}P_{v,\,q}(\Lambda)=
\prod\limits_\phi \left(
M_{v^{deg\,\phi},q^{deg\,\phi}}(\Lambda_\phi)
\prod\limits_{r=1}^\infty\left(1-(vq^{-r})^{deg\,\phi}\right)^{-r}\right).
\end{eqnarray} In particular, as $P_{v,q}$ and $M_{v,q}$ are
probabilistic, one has:
\begin{eqnarray}\label{M u 5}
\prod\limits_{r=0}^\infty\left(1-vq^{-r}\right)^{-1}=
\prod\limits_\phi
\prod\limits_{r=1}^\infty\left(1-(vq^{-r})^{deg\,\phi}\right)^{-r}.
\end{eqnarray}
Now fix $k\in \mathbb{N}$, arbitrary partitions $\lambda_j\in
\mathbb{Y}$ and distinct polynomials $\phi_j\in \Phi,j=1,\ldots,k$.
Using (\ref{M u 3}-\ref{M u 5}) we get:
\begin{eqnarray}\label{M u 6}\begin{split}
   \prod\limits_{i=1}^n(1-q^{-i})^{-1}
   \mu_n(\{\Lambda:\Lambda_{\phi_j}=\lambda_j,j=1,\ldots,k\}\cap
\mathbb{L}_n)=\\
   [v^n]
\prod\limits_{r=0}^\infty\left(1-vq^{-r}\right)^{-1}\prod\limits_{j=1}^k
M_{v^{deg\,\phi_j},q^{deg\,\phi_j}}(\lambda_j).\end{split}
  \end{eqnarray}
To finish the proof of the theorem \ref{main} we need the following
statement (see \cite{F2}, paragraph 2.2, lemma 2):
\begin{Lm}\label{limit}
Let $f(v)$ be a holomorphic in $|v|<1$ function, such that the
Taylor series of $f$ around $0$ converges at $v=1$. Then the next
limit exists:
\begin{eqnarray*}
\lim\limits_{n\rightarrow\infty}[v^n]f(v)(1-v)^{-1}=f(1).
\end{eqnarray*}
\end{Lm}
The proof follows from the identity $(1-v)^{-1}=1+v+v^2+\ldots$ for
$|v|<1$. Using the lemma (\ref{limit}) for \begin{eqnarray*}f(v)=
\prod\limits_{r=1}^\infty\left(1-vq^{-r}\right)^{-1}\prod\limits_{j=1}^k
M_{v^{deg\,\phi_j},q^{deg\,\phi_j}}(\lambda_j)
 \end{eqnarray*} one obtains from (\ref{M u 6}):
 \begin{eqnarray*}
  \lim\limits_{n\rightarrow \infty}
   \mu_n(\{\Lambda:\Lambda_{\phi_j}=\lambda_j,j=1,\ldots,k\}\cap
\mathbb{L}_n)=\prod\limits_{j=1}^k
M_{1,\,q^{deg\,\phi_j}}(\lambda_j).
  \end{eqnarray*}

\medskip
Mathematics Division, Institute for Low Temperature Physics and
Engineering, 47 Lenin ave., Kharkov 61103, Ukraine

e-mail: artemdudko@rambler.ru


\begin{thebibliography}{}

\bibitem{A} Andrews, The theory of partitions. Encyclopedia of Mathematics
and its applications, Vol. 2. Addison-Wesley Publishing Co.,
Reading, Mass.-London-Amsterdam, 1976.

\bibitem{Bor-Ok-Ol} Borodin A., Okounkov A., Olshanski G., {\it Asymptotics of
Plancherel measures for symmetric groups}, J. Amer. Math. Soc. 13
(2000), no. 3, 481--515, math.CO/9905032.

 \bibitem{F1} Fulman J., {\it Probability in the classical groups over
finite fields: symmetric functions, stohastic algorithms and cycle
indeces}, Ph.D. Thesis, Harvard University, 1997.

\bibitem{F2} Fulman J., {\it Random matrix theory over finite
fields}, arXiv:math/0003195.

\bibitem{F3} Fulman J., {\it Convergence rates of random walk on
    irreducible representations of finite groups},
J. Theoret. Probab. Vol. 21, March 2008, arXiv:math/0607399.

\bibitem{Mac} Macdonald I., Symmetric functions
 and Hall Polinomials. Clarendon Press, Oxford. 1995.

\bibitem{Ok}  Okounkov A., {\it Infinite wedge and random partitions},
Selecta Math., New Ser., 7 (2001), 57–81, math.RT/9907127.



\end{thebibliography}
\end{document}